\documentclass[12pt,leqno]{article}  
\usepackage{a4wide,epsf,latexsym}
\usepackage[centertags]{amsmath}
\usepackage{graphicx}
\usepackage{cite}
\numberwithin{equation}{section}

\newtheorem{remark}{Remark}

\newcommand{\ve}{\varepsilon}

\begin{document}

\title {Layer-adapted meshes: Milestones in 50 years of history
}

\author{Hans-G. Roos, TU Dresden}

\date{ August 2019}

\maketitle

\begin{abstract}
50 years ago the first paper on layer-adapted meshes appeared. We sketch the development
in all these years with special emphasis on important ideas.
\end{abstract}


\section{Introduction}
We mainly present the construction of meshes suitable for a one-dimensional
second order convection-diffusion problem posed in $[0,1]$ with a layer term
$E=\exp(-\gamma\,x/\varepsilon)$ (of course, the layer could also be located at $x=1$ or at both
endpoints of the interval, for instance, in case of a reaction-diffusion problem). For 2D problems with boundary layers, tensor product ideas allow
to extend the 1D principles of the mesh construction to the two-dimensional case.

For {\it finite difference methods} it seems natural to reach for \emph{uniform
convergence}\index{uniformly convergent scheme}  in the discrete maximum norm; that is, the
computed solution $\{u_i^N\}_{i=0}^N$ satisfies
\begin{equation}\label{eq1:unifcgce}
\|u-u^N\|_{\infty,d} := \max_{i=0,\dots,N} |u_i-u_i^N| \le C N^{-\alpha}
\end{equation}
for some positive constants $C$ and  $\alpha$ that are independent of $\varepsilon$ and  $N$.
A~power of $N$ is a suitable measure of the error $u-u^N$ for the particular families of meshes
that we will discuss, but a bound of this type is inappropriate for an
arbitrary family of meshes; see \cite{SW96}.

The aim to achieve uniform convergence in the maximum norm is demanding and leads even in 1D to meshes
which sometimes do not have desirable properties. Therefore, we will take into account as well
meshes where the constant $C$ in \eqref{eq1:unifcgce} will weakly depend on $\varepsilon$. Such a
desirable property is, for example, the local quasi-uniformity of the
mesh.

For {\it finite element methods} in 2D it is unrealistic to hope for uniform convergence in the
maximum norm on very non-regular meshes. Therefore, one analyzes finite element methods in
scaled Sobolev space norms. To achieve the $\ve$-independence of the constants in the estimates
again layer-adapted meshes are necessary.

\section{Bakhvalov 1969}
Bakhvalov's original mesh \cite{Ba69}\index{Bakhvalov mesh} uses mesh points $\{x_i\}$
near $x=0$  defined by
$$
    q(1-\exp(-\frac{\gamma x_i}{\mu\varepsilon}))=\frac{i}{N}
$$
 in some interval $[0,\tau]$, here $q$ is a parameter. More precisely, Bakhvalov proposed the mesh-generating function
$$
\phi(t)=\left\{ \begin{array}{ll}
    -\frac{\mu\varepsilon}{\gamma}\ln \frac{q-t}{t}\quad {\rm for}\quad t\in[0,\tau]\\
    \phi(\tau)+\phi'(\tau)(t-\tau)\quad {\rm for}\quad t\in[\tau,1].
    \end{array}\right.
$$
Remarkably,  $\tau$ is defined by the requirement that the mesh-generating function is $C^1$.
Thus, $\tau$ has to solve the nonlinear equation
$$
   \phi'(\tau)=\frac{1-\phi(\tau)}{1-\tau}.
$$
On that mesh, Bakhvalov studied the finite difference method for reaction-diffusion problems
and proved uniform convergence of the second order.

Surprisingly, more than 10 years nobody cited Bakhvalov's paper. In Russian, Vassiljeva was the first to
refer to Bakhvalov in 1982, in English written papers Blatov, Boglaev and Liseikin quoted Bakhvalov 1990 ( papers
of Vulanovic since 1983 based on Bakhvalov's work appeared in a small unknown journal in Novi Sad, Serbia).

Boglaev, Liseikin and Vulanovic presented different versions of Bakhvalov's mesh. For instance,
 a
{\it Bakhvalov-type mesh}\index{ Bakhvalov-type mesh} is given by
\begin{equation}\label{B-type-mesh}
 x_i= -\frac{\mu\varepsilon}{\gamma}\ln\left(1-2(1-\varepsilon)\frac{i}{N}\right),\quad i=0,1,\cdots, N/2.
\end{equation}
In $[\sigma^*,1]$ the mesh is equidistant, where
the transition point from the fine to the coarse mesh is defined by
\begin{equation}\label{transition point B}
\sigma^*=\min\{1/2,\mu \frac{\varepsilon}{\gamma}\ln \frac{1}{\varepsilon}\}.
\end{equation}

Bakhvalov-type meshes are simpler than the Bakhvalov meshes and the mesh-generating function is not
longer $C^1$. But both meshes are not locally quasi-equidistant. {\it In some cases for these meshes and finite difference
methods in 1D
optimal error estimates are known, but the analysis is often more complicated than for
Shishkin-type meshes (see Section 5). For linear finite elements in 1D an optimal error estimate
was proved 2006 in \cite{Ro06}, but is open in 2D (see \cite{RS12})}.

As Lin{\ss} pointed out in \cite{Li10}, a Bakhvalov-mesh can also be generated by equidistributing
the monitor function
$$
M(s)=\max\left(1, \tilde K\gamma\varepsilon^{-1}e^{-\frac{\gamma s}{\mu\varepsilon}}\right).
$$

In several papers (see \cite{LP89} and its references)  Liseikin examines the convergence of
finite difference methods when using mesh generating functions $\lambda(t)$ of the given
independent variable that satisfy $|\lambda'(t)| \le C$ for all $t\in [0,1]$. This approach
generates a graded grid of Bakhvalov type\index{Bakhvalov-type mesh}. His book \cite{Lis99}
develops a general theory of grid generation. The analysis in these sources is written in terms of
``layer-resolving transformations"; their relationship to mesh generating functions in a singular
perturbation context is discussed in \cite{Vul07}.

In \cite{SY12} we find the proposal to generate a mesh by the implicitly defined function
\begin{equation}\label{Lambert}
 \xi(t)-e^{\frac{\gamma\xi(t)}{\mu\varepsilon}}+1-2t=0.
\end{equation}
The mesh has the advantage that it is not necessary to use different mesh generating functions in
different regions. But \eqref{Lambert} is not so easy to solve, however, a solution based on the
use of Lambert's W-function is possible. Some difference schemes (and finite elements) on that
mesh can be analyzed similarly as on Bakhvalov-type meshes \cite{RTU14}.

\section{van Veldhuizen 1978}
Van Veldhuizen \cite{Ve78} studied a second order 1D convection-diffusion problem with constant convection term and
proved for the finite element method with $k$-th order polynomials and Radau quadrature for $k\ge 2$
\begin{equation}
\|u-u^N\|_\infty\le C(h^{k+1}+h_{\ve}^{k+1}+e^{-\delta/\ve}).
\end{equation}
Here the interval $[0,1]$ is decomposed into $[0,\delta]$ and $[\delta,1]$, and $h_{\ve}$ and
$h$ are the maximal step sizes in the corresponding intervals (for $k=1$ he proved a first
order result).

This estimate motivated Veldhuizen to choose
\begin{equation}
\delta=2\ve(2+(k+1)\ln N).
\end{equation}
{\it That means, 10 years before Shishkin Veldhuizen proposed to choose ,,Shishkin's'' transition
point and provided us with the tools to analyse the finite element method on a
Shishkin mesh (see Section 4) }!

Although published in the journal ,,Numerische Mathematik'', almost nobody observed the important
paper of Veldhuizen, and up to today we have only 5 citations. Moreover, in his numerical experiments
Veldhuizen used a {\it Bakhvalov-Shishkin mesh} (we will discuss these meshes in Section 6), without
knowing Bakhvalov's paper from 1969.

\section{Gartland 1988}
Gartland \cite{Ga88} studied higher order finite difference methods in 1D  and graded a mesh in the following way:
$$
x_0=0,\quad x_1=\varepsilon H,\quad x_{i+1}=x_i+h_i
$$
with
\begin{equation}\label{Gartland}
  h_i=\min\left(H, \varepsilon H e^{\frac{\gamma x_i}{2\varepsilon}},e h_{i-1}.
   \right)
\end{equation}
The restriction $h_i\le e h_{i-1}$ ensures that the mesh is {\it locally quasi-equidistant}.

\begin{remark}\label{rem1:upwindarbmesh}
If simple upwinding \index{simple upwind scheme} for a convection-diffusion problem is uniformly
convergent \index{uniformly convergent scheme} in the sense of \eqref{eq1:unifcgce} for some
constant $\alpha >0$, and the mesh is locally quasi-equidistant \index{locally quasi-equidistant
grid} (uniformly in $\varepsilon$), then the number $N$ of mesh intervals must increase as
$\varepsilon\rightarrow 0$. To see this, observe that the arguments of \cite{St03} are still valid
when slightly modified by considering a limit as $N\rightarrow\infty$ with $\varepsilon \ge h_1$
and $i=1$; one then arrives at the conclusion of that paper that $h_1 = o(\varepsilon)$. (There
are some minor extra mesh assumptions such as existence of $\lim_{N\rightarrow\infty} h_1/h_2$ and
$\lim_{N\rightarrow\infty} h_2/h_1$.) But the mesh diameter is at least $1/N$, so the locally
quasi-equidistant property implies that $\varepsilon K^N \ge 1/N$, where $K$ is the constant in
$$
   h_i\le Kh_j \quad {\rm for}\quad |i-j|\le 1.
$$
 Hence $NK^N \ge 1/\varepsilon$, so $N \approx \log_K(1/\varepsilon)$.
 \hfill$\clubsuit$
\end{remark}

Introducing the transition points $x^*,x'$ by
$$
    x^*\approx K\varepsilon \ln \frac{K}{H}, \quad x'\approx K\varepsilon \ln \frac{K}{\varepsilon}
$$
Gartland observed that the number of mesh points in the inner region $[0,x^*]$ as well in the
outer region $[x',1]$ is of order $O(1/H)$, but in the transition region $[x^*,x']$ of order
$O(\ln\ln \frac{H}{\varepsilon})$.
On that mesh Gartland proved the uniform convergence of some finite difference schemes.

We call the modification of the mesh where \eqref{Gartland} is replaced by
\begin{equation}
 h_i=\min\left(H, \varepsilon H e^{\frac{\gamma x_i}{2\varepsilon}}.
   \right)
   \end{equation}
{\it Gartland-type mesh}\index{Gartland-type mesh}. The number of mesh points is now independent
of $\varepsilon$ and the mesh allows {\it optimal error estimates}. The mesh is not locally
quasi-eqidistant.

Finite element methods on
 Gartland-type meshes were first studied 1997
in \cite{RS97}, see also \cite{RTU15}. There is also a close relation to the results of Liu and Xu
\cite{LX06,LX09,LWX16}.

\section{Shishkin 1988}
Shishkin spread the great idea to use the very simple {\it piecewise constant} meshes
in combination with
the transition point $\sigma$ from the fine to the coarse mesh
 defined by
\begin{equation}\label{transition point S}
\sigma=\min\{1/2,\mu \frac{\varepsilon}{\gamma}\ln  N \}.
\end{equation}
Consequently, for small $\varepsilon$, one has $E(\sigma)=N^{-\mu}$, and $\mu$ is
chosen in dependence of the order of the method used.

Then, each of the intervals $[0,\sigma]$ and $[\sigma,1]$ is subdivided  equidistantly  in
$N/2$ subintervals.
It is not vital  that one has exactly the same number of
subintervals in $[0,\sigma]$ and $[\sigma,1]$. All that the
theory demands is that as $N\rightarrow\infty$ the number of
subintervals in each of these two intervals is bounded below by
$CN$ for some constant $C>0$.

The coarse  part of this Shishkin mesh has spacing
$H=2(1-\sigma)/N$, so $N^{-1}\le H \le 2N^{-1}$. The fine part has
spacing $h=2\sigma/N = (4/\gamma)\varepsilon N^{-1}\ln N$, so
$h\ll\varepsilon$. {\it Thus there is a very abrupt change in mesh size
as one passes from the coarse part to the fine part. The mesh is
not locally quasi-equidistant,  uniformly in $\varepsilon$}.

Before Shishkin's work uniform finite difference methods were often based on the properties
of pointwise uniform consistency and uniform stability in the maximum norm. But pointwise
uniform consistency is not necessary. While Shishkin used the maximum principle and barrier functions
in the analysis, other approaches use improved stability properties (Andreev, Kopteva, Lin{\ss} 1996-98
\cite{AK96,AK98,Li10}).

In 1D, we know uniform second order finite difference schemes on Shishkin meshes, but in 2D
mostly only the first order upwind scheme is analyzed (for convection diffusion problems).

Finite element methods on Shishkin meshes in 1D were first studied 1995 by Sun and Stynes  \cite{SS95b}, the
analysis for second order problems was also published in the two books  \cite{MOS96,RST96} from 1996. But
the more important analysis in 2D was still in open problem, because in 2D the analysis of finite element
methods is traditionally based on shape regular or isotropic meshes (see Section 6).

If a method for a problem with a smooth solution has the order $\alpha$, due to the fine mesh
size $h=O(\varepsilon N^{-1}\ln N)$ in the case $u^{(k)}\approx \varepsilon ^{-k}$ we can expect
 that the error on a Shishkin mesh is of the order $O((N^{-1}\ln N)^\alpha)$. Especially for higher
order methods the logarithmic factor is troublesome. An optimal mesh should generate an error of the
order $O(N^{-\alpha})$. This is the reason for introducing S-type meshes, see Section 7.

Shishkin meshes and S-type meshes have many advantages, but there are also disadvantages: the
robustness (see Section 9) and the loss of the possibility to use some typical ingredients
in the analysis of finite element methods on isotropic or local uniform meshes.

\section{Apel and Dobrowolski 1992/97: anisotropic meshes}
In the classical theory of finite element methods one uses the Lagrange interpolant
$u^I\in P_k$ on some element $K$ for $u\in W^{k+1,p}(K)$ and its approximation
property
\[
|u-u^I|_{m,p,K}\le C\frac{h_K^{k+1}}{\rho_K^m}|u|_{k+1,p,K}.
\]
Here $h_K$ is the diameter of $K$ and $\rho_K$ the length of the largest ball
inscribed $K$. The classical theory {\it assumes} a bounded aspect ratio
\[
\frac{h_K}{\rho_K} \le C,
\]
which excludes anisotropic elements.

In 1992 Apel and Dobrowolski \cite{AD92} proved sharp anisotropic interpolation
error estimates in the case $m=1$, Apel and Lube extended 1994 the results to
general $m$. Finally in 1999 Apel presented a general theory of anisotropic
elements in his famous book \cite{Ape99}.

We just sketch in a simple situation an anisotropic result:
Suppose that each element $K$ (triangle or rectangle) of a mesh is
contained in a rectangle with side lengths $(h_x,h_y)$ and
contains a rectangle with side lengths $(C_2h_x,C_2h_y)$ for some
fixed constant $C_2>0$. In the case of triangles, assume also a
maximum angle condition: the interior angles of every mesh
triangle are bounded away from $\pi$. (Triangular \index{Shishkin
mesh}Shishkin meshes have maximum angle $\pi/2$ and consequently
satisfy this condition.) Then there exists a constant $C$ such
that
\begin{subequations}\label{eq3:Apel}
    \begin{align}
\|v-v^I\|_{0,p,K} &\le  C \sum_{|\alpha|= m}h^{\alpha} \|D^\alpha v\|_{0,p,K}
    \quad\text{for } m=1,2,    \label{eq3:Apela}\\
\|\partial_x(v-v^I)\|_{0,p,K} &\le
  C \sum_{|\alpha| =1}h^{\alpha} \|D^\alpha \partial_x v\|_{0,p,K},
    \label{eq3:Apelb}\\
\|\partial_y(v-v^I)\|_{0,p,K} &\le
  C \sum_{|\alpha| =1}h^{\alpha} \|D^\alpha \partial_y
  v\|_{0,p,K}, \label{eq3:Apelc}
    \end{align}
\end{subequations}
where we set $\ h^{\alpha} = h_x^{\alpha_1} h_y^{\alpha_2}$ and $\ D^{\alpha} = {\partial_x}^{\alpha_1} {\partial_y}^{\alpha_2}$ .

{\it 1997 Dobrowolski created the idea to use these estimates for the first analysis of a
singularly perturbed convection-diffusion problem, using linear or bilinear finite
elements on a Shishkin mesh in 2D} \cite{DR97}. Parallel and independently Stynes and O'Riordan used
a very special technique for bilinear elements \cite{SO97}.

Since that time anisotropic interpolation error estimates are a standard ingredient
to analyze finite element methods on layer-adapted meshes.

\section{ S-type meshes 1999}
In 1999 Lin{\ss} introduced Bakhvalov-Shishkin meshes similarly as van Veldhuizen 20 years before
and analyzed  finite difference and finite element methods. More generally, we defined
in \cite{RL99} the class of S-type meshes. In the fine subinterval $[0,\sigma]$ with the transition point
$\sigma=\frac{\mu\varepsilon}{\gamma}\ln N$ from the fine to the coarse mesh
we use a {\it mesh-generating function}. Assuming the
function
$\lambda:\,\, [0,1/2]\mapsto [0, \ln N]$ to be strictly increasing, set
$$
   x_i= \frac{\mu\varepsilon}{\gamma}\lambda(i/N),\quad i=0,1,\cdots, N/2.
$$
We call such meshes {\it Shishkin-type meshes}\index{Shishkin-type mesh}.
It turns out  that in error estimates for Shishkin-type meshes often the factor $\max|\psi'(\cdot)|$
appears, where $\psi$ is the {\it mesh-characterizing function} defined by
$$
   \psi:=e^{-\lambda}:\quad [0,1/2]\mapsto [1, 1/N].
$$
For the original Shishkin mesh we have $\max|\psi'(\cdot)|=O(\ln N)$. A popular optimal mesh
is the {\it Bakhvalov-Shishkin mesh}\index{Bakhvalov-Shishkin mesh} with
$$
 \psi(t)=1-2t(1-N^{-1}) \quad {\rm and}\quad \max|\psi'(\cdot)|\le 2.
$$
The mesh points of the fine mesh are given by
\begin{equation}\label{BS-mesh}
 x_i= -\frac{\mu\varepsilon}{\gamma}\ln\left(1-2(1-N^{-1})\frac{i}{N}\right),\quad i=0,1,\cdots, N/2.
\end{equation}
Another optimal mesh is the {\it Vulanovic-Shishkin} mesh, for details
and other possibilities to choose $\lambda$, see \cite{Li10}.

For higher order finite elements one gets for optimal meshes with bounded $\max|\psi'(\cdot)|$ the
optimal error estimate in the energy norm
\begin{equation}
\|u-u^N\|_\ve \le CN^{-k},
\end{equation}
which is for higher $k$ much better than the result on a Shishkin mesh. One can also try to optimize
$\psi$ or the error constants in the estimates \cite{RTU15}.

Further modifications of Shishkin meshes due to Vulanovic are also described in \cite{Li10}.
In \cite{FX17} from 2017 we find  a generalization of Shishkin-type meshes based on the property\\
$\lambda(1/2)=\ln (\theta N)$ with some additional parameter $\theta$. This allows to
characterize the so called eXp-mesh \cite{Xe02} from Xenophontos 2002 as generalized Shishkin-type mesh.

\section{Kopteva and Stynes 2001: Adaptively generated meshes}
For a long time {\it monitor functions} are  a standard tool to generate meshes, for instance,
the function
\[
   M=\sqrt{1+(u')^2}
\]
related to the arc length. Several authors studied adaptive algorithms for singularly
perturbed problems based on that monitor function \cite{HRR94,Mac99,BM00,QS99}.

The breakthrough came with the results 2001 of Kopteva on {\it a posteriori error bounds} for
some numerical methods for convection-diffusion problems in 1D. For a
conservative form of the upwind finite difference method Kopteva proved
\[
 \|u^N-u\|_\infty\le C\max_i h_i \sqrt{1+ (D^-u_i^N)^2}\,.
\]
Here $u^N$ is the linear interpolant of the computed solution.

Setting $M_i := \sqrt{1+ (D^-u_i^N)^2}$ for $i=1,\dots, N$, Stynes and Kopteva introduced
 the  \emph{equidistribution problem}: Find $\{ (x_i, u_i^N)\}$, with
the $\{u_i^N\}$ computed from the $\{x_i\}$ by means of the upwind scheme,  such
that
\begin{equation}\label{eq:equi}
h_i M_i = \frac{1}{N}\sum_{j=1}^N h_j M_j  \ {\rm for\
}i=1,2,\ldots,N.
\end{equation}

Unfortunately, this is a nonlinear problem. But Stynes and Kopteva presented an algorithm
which stops in less than
 $ K \le C_3 |\ln \varepsilon|/(\ln N)$ steps, such that
 \[
\|e^{(K)}\|_\infty \le C_4N^{-1},
\]
where $e^{(k)}$ is the error in the $k^{\rm th}$ solution
computed by the algorithm.

Experimental evidence shows that the final mesh computed by the algorithm  is strikingly close to
a Bakhvalov mesh\index{Bakhvalov mesh} inside the boundary layer; see \cite[Fig.~2]{KS01b}. In
contrast, most adaptive algorithms will not generate a mesh resembling a
Shishkin mesh.

{\it Unfortunately, in 2D the situation is very different}.

\noindent An adaptive procedure designed for problems with layers should include
an anisotropic refinement strategy. While several  anisotropic mesh
adaptation strategies do exist, all are more or less
heuristic. {\it We do not know of  any strategy for convection-diffusion problems in two dimensions
where it is proved that, starting from some standard mesh, the refinement strategy is guaranteed
to lead to a mesh that allows robust error estimates}.

A necessary tool is a robust error estimator on anisotropic meshes; we shortly sketch the situation.
Very important is the {\it robustness} of the estimators, even on isotropic meshes.
In \cite{San01} Sangalli proves the robustness of a certain  estimator
for the residual-free bubble method  applied to convection-diffusion
problems. The analysis uses the norm
\begin{equation}\label{eq3:3.533o}
||w||_{San}:=\|w\|_\varepsilon +\|b\cdot \nabla w\|_*,\quad\text{where}\quad \|\varphi\|_*=\sup
\frac{\langle\varphi,v\rangle}{\|v\|_\varepsilon}\,.
\end{equation}
Although Sangalli's approach is devoted to residual-free bubbles, the same analysis works for the
Galerkin method and the SDFEM. For the convection-diffusion problem, the residual error estimator
 is  robust with respect to the dual norm; see \cite{Ve04}.

 The norm $||\cdot||_{San}$ above is defined only implicitly by an
infinite-dimensional variational problem and cannot be computed exactly in practice.
In \cite{San08}  Sangalli pointed out that the norm (\ref{eq3:3.533o}) seems to be  not optimal in
the convection-dominated regime. He proposes
an improved estimator that is robust with respect to this natural norm \cite{San05} for the
advection-diffusion operator,  but studied only the one-dimensional case.
 The relation to another new improved dual norm is studied in detail in \cite{DZ15}.

 Today the dual norm or its modification plays an important role in many papers on
robust a posteriori error estimation for convection-diffusion problems.

In \cite{TV15}, Tobiska and Verf{\"u}rth proved in the dual norm that the same robust a posteriori
error estimator can be used for a range of stabilized methods such as streamline diffusion, local
projection schemes, subgrid-scale techniques and continuous interior penalty methods. Nonconforming
methods are studied in \cite{ZC14}. Variants of discontinuous Galerkin methods are discussed in
\cite{ESV10,GSZ14,LPP16,SZ09,ZS11}. Vohralik \cite{Vo12} presents a very general concept of
a posteriori error estimation based on potential and flux reconstructions.

{\it Most papers mentioned assume isotropic meshes, but
Kopteva designed starting in 2015 different estimators ( residual \cite{Ko15,Ko17}, flux equilibration \cite{Ko18})
for {\it anisotropic meshes}}.

For reaction-diffusion problems
it is unclear that the energy norm is a
suitable norm for these problems because for small~$\varepsilon$ it is unable to distinguish
between the typical layer function of reaction-diffusion problems and zero. {\it It would be desirable
to get robust \emph{a posteriori} error estimates in a stronger norm,  for instance, some
balanced norm or the $L_\infty$ norm}.

The  first result with respect to the maximum norm is the \emph{a posteriori} error estimate of Kopteva
\cite{Ko07b} in 2008 for the standard finite difference method on an arbitrary
rectangular mesh.  Next we sketch the ideas of \cite{DK16} from Demlov and Kopteva for a posteriori error estimation for finite elements
of arbitrary order on {\it isotropic} meshes in the maximum norm 2016.

Using the Green's function of the continuous operator with respect to a point $x$, the error in that point
can be represented by
\[
e(x)=\varepsilon^2(\nabla u_h,\nabla G)+(cu_h,G).
\]
For some $G_h\in V_h$ we obtain
\[
e(x)=\varepsilon^2(\nabla u_h,\nabla (G-G_h))+(cu_h,G-G_h).
\]
Integration by parts yields
\[
e(x)=\frac{1}{2}\sum_{T\in {\cal T}_h}\int_{\partial T}\varepsilon^2(G-G_h)n_T\cdot[\nabla u_h]
     +\sum_{T\in {\cal T}_h}(cu_h-f-\varepsilon^2\triangle u_h,G-G_h)_T.
\]
Choosing for $G_h$ the Scott-Zhang interpolant
 of $G$, one needs sharp estimates for $G$ to control
the interpolation error. These are
 collected in Theorem 1 of \cite{DK16}. Thus, one obtains finally with
 $l_{h}:=\ln(2+\tilde \varepsilon {\underline h}_{}^{-1} )$ (the constant $\tilde \varepsilon$ is
 of order $\varepsilon$ and $\underline h=\min h_T$)
 \begin{subequations}\label{Kop1}
\begin{align}
\|u-u_h\|_\infty &\le C\max_{T\in {\cal T}_h}(
   \min(\tilde\varepsilon,l_hh_T)\|[\nabla u_h]\|_{\infty,\partial T}\\
  &+ \min(1,l_hh_T^2\varepsilon^{-2})\|cu_h-f-\varepsilon^2\triangle u_h\|_{\infty, T})
.\nonumber
\end{align}
\end{subequations}

On {\it anisotropic meshes}, \index{anisotropic mesh} in 2015 Kopteva  also derived an a posteriori error estimator in the
maximum norm \cite{Ko15}, now for {\it linear} finite elements. Suppose that the
triangulation satisfies the maximum angle condition. Then the first result of \cite{Ko15}
gives
\begin{subequations}\label{Kop2}
\begin{align}
\|u-u_h\|_\infty &\le C\,l_h\max_{z\in \cal N}(
   \min(\varepsilon,h_z)\|[\nabla u_h]\|_{\infty,\partial \omega_z}\\
  &+ \min(1,h_z^2\varepsilon^{-2})\|cu_h-f\|_{\infty, \omega_z})
.\nonumber
\end{align}
\end{subequations}
Here $\omega_z$ is the patch of the elements surrounding some knot $z$ of the triangulation,
$h_z$ the diameter of  $\omega_z$. In a further estimator the second term of \eqref{Kop2},
which has isotropic character, is replaced by a sharper result with more anisotropic
nature.

To prove \eqref{Kop2} two difficulties arise. First, it is necessary to use scaled trace
inequalities. Moreover, instead of using the Scott-Zhang interpolant of the Green's function (whose
applicability is restricted on anisotropic meshes) Kopteva uses some standard Lagrange
interpolant for some continuous approximation of $G$. But the construction is based on the
following additional assumption on the mesh. Let us introduce $\Omega_1:=\{T: h_T\ge c_1\varepsilon\}$
and  $\Omega_2:=\{T: h_T\le c_2\varepsilon\}$ with some positive $c_1<c_2$. Then, the additional
assumption requires that the distance of $\Omega_1$ and $\Omega_2$ is at least some $c_3\varepsilon$
with $c_3>0$.

The last condition excludes an too abrupt change of the mesh size, typically for Shishkin meshes.
But other layer-adapted meshes satisfy that condition, for instance, Bakhvalov meshes or
Bakhvalov-Shishkin meshes.

Unfortunately, we still miss an adaptive strategy based on these estimators leading to optimal
meshes.

\section{Duran and Lombardi 2006: simple recursively graded meshes}
Very simple is the {\it Duran-Lombardi mesh} \cite{DL06}\index{Duran-Lombardi mesh} defined by
$$
\begin{array}{ll}
x_0=0,\,\,x_i=i\kappa H\varepsilon\quad {\rm for}\quad 1\le i\le \frac{1}{\kappa H}+1\\
 x_{i+1}=x_i+\kappa H x_i\varepsilon\quad {\rm for}\quad  \frac{1}{\kappa H}+1\le i\le M-2,\quad x_M=1.
\end{array}
$$
Here M is chosen such that $x_{M-1}<1$ but $x_{M-1}+\kappa H x_{M-1}\ge 1$, assuming that the last interval
is not extremely small.

The mesh is locally quasi-equidistant and glitters by its simplicity. Almost uniform error estimates with respect to $H$ are possible, but the number of mesh-points is
proportional to $\frac{1}{H}\ln \frac{1}{\varepsilon}$. The finite element analysis presented in \cite{DL06}
differs a bit from the analysis on a S-type mesh. It has the advantage to require only estimates for
derivatives (in contrast to a solution decomposition with estimates for the components of the decomposition).

For finite difference methods an analysis on a DL mesh is not known. But for the upwind method, using
\[
\|u-u^H\|_{\infty,d}\le C \max\int_{x_{k-1}}^{x_k}(1+|u'|)
\]
(see \cite{Li10})
one can simple conclude
\[
\|u-u^H\|_{\infty,d}\le C\,H.
\]

Practically, the mesh has two important advantages in comparison to S-type meshes.
First, there is no need to define a transition
point. And, remarkably, the mesh has the following robustness property: a mesh defined for some
$\ve^*$ can also be used for larger values of the parameter. More precisely: If $\ve$ is the
perturbation parameter in the equation which varies, we define a mesh based on a smaller value
$\ve^*$. Then in a certain range the error is smaller than using $\ve$ to define the mesh.
Numerical experiments and some heuristic arguments show optimality for $\ve^*\approx \ve/2$.

Of course, one can also define a mesh  by a general recursive formula as
\[
x_{i+1}=x_i+g(\ve,H,x_i)
\]
and analyze the necessary properties of $g$ to obtain nice convergence results.
One possibility is to define a mesh-generating function by interpolation
of the values given in the mesh points and then to analyze discretization
methods similarly as methods on a S-type mesh, see, for instance, \cite{RTU15}.

\section{What else is there?}
\subsection{ Emelyanov and Sidorov Grids}
Emelyanov used 1995 the ,,optimal'' grids developed by Sidorov 1966 to prove the uniform
convergence of some difference schemes \cite{Em95}. Sidorov's intention was to construct meshes
with the property
\begin{equation}\label{em}
 \sum(\frac{h_{i+1}}{h_i}-1)^2 \Longrightarrow  {\rm min.,}
\end{equation}
Instead solving \eqref{em}, he solved introducing $x_i=x(i)$ the continuous problem
\begin{equation}
 \int_0^N\frac{x_{\xi\xi}^2}{x_\xi^2} \Longrightarrow  {\rm min.}\quad {\rm with}\,
 \int_0^Nx_\xi=1, x_\xi(0)=A, x_\xi(N)=B.
\end{equation}
The exact solution of that problem is known. Moreover, Sidorov meshes have the nice property
$h_i-h_{i-1}=O(N^{-2})$, useful for the analysis of difference schemes.

Emelyanov used a fine subinterval at the layer and a Sidorov mesh there and a coarse subinterval
with $h_i=O(N^{-1})$ and $h_i-h_{i-1}=O(N^{-2})$ to prove uniform convergence. It would be interesting
to know whether or not one can choose in a convection-diffusion problem the first and the last mesh size in such a way that the
direct application of Sidorov's approach leads to uniform convergence.

\subsection{Admissibly graded meshes by Liu and Xu}
Liu and Xu study 2009 the Galerkin method for a one-dimensional 2m-th order convection-diffusion problem with
Hermite splines of degree $2r-1$ \cite{LX09}. For simplicity, we sketch their ideas for the mesh construction
in the case $m=r=1$. The mesh is called {\it admissible}, if $h_i\le CN^{-1}$ and
\begin{equation}\label{con}
\sum(\frac{h_i}{\ve})^3e^{-2x_{i-1}/\ve}\le CN^{-2}.
\end{equation}
For other values of $m,r$ a further condition is required. The condition \eqref{con} is sufficient to prove
that the linear interpolant of the layer part $E$ satisfies
\[
\ve^{1/2}|E-E^I|_1\le CN^{-1},\quad \ve^{-1/2}\|E-E^I\|_0\le CN^{-1}.
\]
Next, Liu and Xu show that the condition
\begin{equation}\label{con2}
h_i\le \min\left(S\ve N^{-1}e^{x_{i-1}/(2\ve)},N^{-1} \right)\quad {\rm for}\,\,i\in {\cal N}_{N'}
\end{equation}
with $N'\le CN$ is sufficient for the admissibility of the mesh. Condition \eqref{con2} remembers us
of a Gartland-type mesh.

2016 in \cite{LWX16} Li, Wu and Xu (for a two-dimensional, second order reaction-diffusion problem)
present an explicit realization of a mesh satisfying \eqref{con2}. The result is a Bakhvalov mesh.

\subsection{hp meshes}
When analyzing hp finite element methods for singularly perturbed problems it is common to use
an hp boundary layer mesh, see \cite{Me02}. For such methods it is possible to prove
exponential convergence.


\begin{thebibliography}{10}

\bibitem{AK96}
V.~B. Andreev and N.~V. Kopteva.
\newblock Investigation of difference schemes with an approximation of the
  first derivative by a central difference relation.
\newblock {\em Comput. Math. Math. Phys.}, 36:1065--1078, 1996.

\bibitem{AK98}
V.~B. Andreev and N.~V. Kopteva.
\newblock On the convergence, uniform with respect to a small parameter, of
  monotone three-point difference schemes.
\newblock {\em Differential Equations}, 34:921--929, 1998.

\bibitem{Ape99}
T.~Apel.
\newblock {\em Anisotropic finite elements}.
\newblock Wiley-Teubner, Stuttgart, 1999.

\bibitem{AD92}
T.~Apel and M.~Dobrowolski.
\newblock Anisotropic interpolation with applications to the finite element
  method.
\newblock {\em Computing}, 47:277--293, 1992.

\bibitem{Ba69}
A.~S. Bakhvalov.
\newblock On the optimization of methods for solving boundary value problems
  with boundary layers (in {R}ussian).
\newblock {\em Zh. Vychisl. Mat. i Mat. Fis.}, 9:841--859, 1969.

\bibitem{BM00}
G.~Beckett and J.~A. Mackenzie.
\newblock Convergence analysis of finite difference approximations on
  equidistributed grids to a singularly perturbed boundary value problem.
\newblock {\em Appl. Numer. Math.}, 35(2):87--109, 2000.

\bibitem{DK16}
A.~Demlov and N.~Kopteva.
\newblock Maximum-norm a posteriori error estimates for singularly perturbed
  elliptic reaction-diffusion problems.
\newblock {\em Numer. Math.}, 133:707--742, 2016.

\bibitem{DR97}
M.~Dobrowolski and H.-G. Roos.
\newblock A priori estimates for the solution of convection-diffusion problems
  and interpolation on {S}hishkin meshes.
\newblock {\em Z. Anal. Anwendungen}, 16:1001--1012, 1997.

\bibitem{DZ15}
S.~Du and Z.~Zhang.
\newblock A robust residual-type a posteriori error estimator for
  convection-diffusion equations.
\newblock {\em J. Sci. Comput.}, 65:138--170, 2015.

\bibitem{DL06}
R.~G. Dur{\'a}n and A.~L. Lombardi.
\newblock Finite element approximation of convection diffusion problems using
  graded meshes.
\newblock {\em Appl. Numer. Math.}, 56:1314--1325, 2006.

\bibitem{Em95}
K.~V. Emelyanov.
\newblock On optimal grids and their application to the solution of problems
  with a singular perturbation.
\newblock {\em Russ. J. Num. Anal. Math. Modelling}, 10:299--309, 1995.

\bibitem{ESV10}
A.~Ern, A.~F. Stephansen, and M.~Vohralik.
\newblock Guaranteed and robust discontinuous {G}alerkin a posteriori error
  estimates for convection-diffusion-reaction problems.
\newblock {\em J. Comput. Appl. Math.}, 234:114--130, 2010.

\bibitem{FX17}
S.~Franz and C.~Xenophontos.
\newblock A short note on the connection between layer-adapted exponentially
  graded and {S}-type meshes.
\newblock {\em CMAM}, 18:199--202, 2017.

\bibitem{Ga88}
E.~C. Gartland.
\newblock Graded-mesh difference schemes for singularly perturbed two-point
  boundary value problems.
\newblock {\em Math. Comp.}, 51:631--657, 1988.

\bibitem{GSZ14}
S.~Giani, D.~Sch\"otzau, and L.~Zhu.
\newblock An a posteriori error estimate for hg-adaptive {D}{G} methods for
  convection-diffusion problems on anisotropically refined meshes.
\newblock {\em Computers and Math. with Appl.}, 67:869--887, 2014.

\bibitem{HRR94}
W.~Huang, Y.~Ren, and R.~D. Russell.
\newblock Moving mesh partial differential equations (mmpdes) based on the
  equidistribution principle.
\newblock {\em SIAM J. Numer. Anal.}, 31:709--731, 1994.

\bibitem{Ko07b}
N.~Kopteva.
\newblock Maximum norm a posteriori error estimate for a 2d singularly
  perturbed reaction-diffusion problem.
\newblock {\em SIAM J. Numer. Anal.}, 2008.
\newblock Electronic publication April 11.

\bibitem{Ko15}
N.~Kopteva.
\newblock Maximum-norm a posteriori error estimates for singularly perturbed
  reaction-diffusion problems on anisotropic meshes.
\newblock {\em SIAM J. Num. Anal.}, 53:2519--2544, 2015.

\bibitem{Ko17}
N.~Kopteva.
\newblock Energy norm a posteriori error estimates for singularly perturbed
  reaction-diffusion problems on anisotropic meshes.
\newblock {\em Numer. Math.}, 137:607--642, 2017.

\bibitem{Ko18}
N.~Kopteva.
\newblock Fully computable a posteriori error estimator using anisotropic flux
  equilibration on anisotropic meshes.
\newblock {\em to appear}, pages~--, 2018.

\bibitem{KS01b}
N.~Kopteva and M.~Stynes.
\newblock A robust adaptive method for a quasilinear one-dimensional
  convection-diffusion problem.
\newblock {\em SIAM J. Numer. Anal.}, 39:1446--1467, 2001.

\bibitem{Li10}
T.~Lin\ss.
\newblock {\em Layer-adapted meshes for reaction-convection-diffusion
  problems}.
\newblock Springer, Berlin, 2010.

\bibitem{Lis99}
V.~D. Lise{\u\i}kin.
\newblock {\em Grid generation methods}.
\newblock Scientific Computation. Springer-Verlag, Berlin, 1999.

\bibitem{LP89}
W.~D. Lise{\u\i}kin and W.~E. Petrenko.
\newblock An adaptive-invariant method for the numerical solution of problems
  with boundary and interior layers (in {R}ussian).
\newblock Computer Center, Academy of Sci., Novosibirsk, 1989.

\bibitem{LWX16}
S.-T. Liu, B.~Wu, and Y.~Xu.
\newblock High order {G}alerkin methods with graded meshes for two-dimensional
  reaction-diffusion problems.
\newblock {\em Int. J. Num. Anal. Mod.}, 13:319--343, 2016.

\bibitem{LX06}
S.-T. Liu and Y.~Xu.
\newblock Galerkin methods based on {H}ermite splines for singular perturbation
  problems.
\newblock {\em SIAM J. Num. Anal.}, 43:2607--2623, 2006.

\bibitem{LX09}
S.-T. Liu and Y.~Xu.
\newblock Graded {G}alerkin methods for the high-order convection-diffusion
  problem.
\newblock {\em Num. Meth. Part. Diff. Equ.}, 25:1262--1282, 2009.

\bibitem{LPP16}
A.~L. Lombardi, P.~Pietra, and M.~Prieto.
\newblock A posteriori error estimator for exponentially fitted discontinuous
  {G}alerkin approximation of advection dominated problems.
\newblock {\em Calcolo}, 53:83--103, 2016.

\bibitem{Ve78}
van M.~Veldhuizen.
\newblock Higher order methods for a singularly perturbed problem.
\newblock {\em Numer. Math.}, 30:267--279, 1978.

\bibitem{Mac99}
J.~Mackenzie.
\newblock Uniform convergence analysis of an upwind finite-difference
  approximation of a convection-diffusion boundary value problem on an adaptive
  grid.
\newblock {\em IMA J. Numer. Anal.}, 19(2):233--249, 1999.

\bibitem{Me02}
J.~M. Melenk.
\newblock {\em $hp$-Finite Element Methods for Singular Perturbations}.
\newblock Springer, Heidelberg, 2002.

\bibitem{MOS96}
J.~J.~H. Miller, E.~O'Riordan, and G.~I. Shishkin.
\newblock {\em Fitted numerical methods for singular perturbation problems}.
\newblock World Scientific Publishing Co. Inc., River Edge, NJ, 1996.

\bibitem{QS99}
Y.~Qiu and D.~M. Sloan.
\newblock Analysis of difference approximations to a singularly perturbed
  two-point boundary value problem on an adaptively generated grid.
\newblock {\em J. Comput. Appl. Math.}, 101(1-2):1--25, 1999.

\bibitem{Ro06}
H.-G. Roos.
\newblock Error estimates for linear finite elements on bakhvalov type meshes.
\newblock {\em Appl. of Math.}, 51:63--72, 2006.

\bibitem{RL99}
H.-G. Roos and T.~Lin{\ss}.
\newblock Sufficient conditions for uniform convergence on layer-adapted grids.
\newblock {\em Computing}, 63:27--45, 1999.

\bibitem{RS12}
H.-G. Roos and M.~Schopf.
\newblock Analysis of finite element methods on bakhvalov-type meshes for
  linear convection-diffusion problems in 2{D}.
\newblock {\em Appl. of Math.}, 57:97--108, 2012.

\bibitem{RS97}
H.-G. Roos and T.~Skalick\'y.
\newblock A comparison of the finite element method on {S}hishkin and
  {G}artland-type meshes for convection-diffusion problems.
\newblock {\em CWI Quarterly}, 10:277--300, 1997.

\bibitem{RST96}
H.-G. Roos, M.~Stynes, and L.~Tobiska.
\newblock {\em Numerical methods for singularly perturbed differential
  equations}, volume~24 of {\em Springer Series in Computational Mathematics}.
\newblock Springer-Verlag, Berlin, 1996.
\newblock Convection-diffusion and flow problems.

\bibitem{RTU14}
H.-G. Roos, L.~Teofanov, and Z.~Uzelac.
\newblock A modified {B}akhvalov mesh.
\newblock {\em Appl. Math. Letters}, 31:7--11, 2014.

\bibitem{RTU15}
H.-G. Roos, L.~Teofanov, and Z.~Uzelac.
\newblock Graded meshes for higher order {F}{E}{M}.
\newblock {\em J. Comput. Math.}, 33:1--16, 2015.

\bibitem{San01}
G.~Sangalli.
\newblock A robust a posteriori estimate for the residual free bubbles method
  aplied to advection-dominated problems.
\newblock {\em Numer. Math.}, 89:379--399, 2001.

\bibitem{San05}
G.~Sangalli.
\newblock A uniform analysis of nonsymmetric and coercive linear operators.
\newblock {\em SIAM J. Math. Anal.}, 36(6):2033--2048, 2005.

\bibitem{San08}
G.~Sangalli.
\newblock Robust a-posteriori estimator for advection-diffusion-reaction
  problems.
\newblock {\em Math. Comp.}, 77(261):41--70 (electronic), 2008.

\bibitem{SZ09}
D.~Sch\"otzau and L.~Zhu.
\newblock A robust a posteriori error estimator for discontinuous {G}alerkin
  methods for convection-diffusion equations.
\newblock {\em Appl. Num. Math.}, 59:2236--2255, 2009.

\bibitem{SW96}
P.~M. Selwood and A.~J. Wathen.
\newblock Convergence rates and classification for one-dimensional
  finite-element meshes.
\newblock {\em IMA J. Numer. Anal.}, 16:65--74, 1996.

\bibitem{SY12}
G.~Soederlind and A.~S. Yadaw.
\newblock The impact of smooth {W}-grids in the numerical solution of singular
  perturbation two-point boundary value problems.
\newblock {\em Appl. Math. Comp.}, 218:6045--6055, 2012.

\bibitem{St03}
M.~Stynes.
\newblock A jejune heuristic mesh theorem.
\newblock {\em Comput. Meth. Appl. Math.}, 3:488--492, 2003.

\bibitem{SO97}
M.~Stynes and E.~O'Riordan.
\newblock A uniformly convergent {G}alerkin method on a {S}hishkin mesh for a
  convection-diffusion problem.
\newblock {\em J. Math. Anal. Appl.}, 214:36--54, 1997.

\bibitem{SS95b}
G.~Sun and M.~Stynes.
\newblock Finite element methods for singularly perturbed higher order elliptic
  two-point boundary value problems {II}:~convection-diffusion type.
\newblock {\em IMA J. Numer. Anal.}, 15:197--219, 1995.

\bibitem{TV15}
L.~Tobiska and R.~Verf\"urth.
\newblock Robust a posteriori error estimates for stabilized finite element
  methods.
\newblock {\em IMA J. Num. Anal.}, 35:1652--1671, 2015.

\bibitem{Ve04}
R.~Verf\"urth.
\newblock Robust a posteriori error estimates for stationary
  convection-diffusion equations.
\newblock {\em SIAM J. Numer. Anal.}, 43:1766--1782, 2005.

\bibitem{Vo12}
M.~Vohralik.
\newblock A posteriori error estimates for efficiency and error control in
  numerical simulations.
\newblock Lecture Notes, University Pierre et Marie Curie, 2012.

\bibitem{Vul07}
R.~Vulanovi{\'c}.
\newblock The layer-resolving transformation and mesh generation for
  quasilinear singular perturbation problems.
\newblock {\em J. Comput. Appl. Math.}, 203:177--189, 2007.

\bibitem{Xe02}
C.~Xenophontos.
\newblock Optimal mesh design for the finite element approximation of
  reaction-diffusion problems.
\newblock {\em IJNME}, 53:929--943, 2002.

\bibitem{ZC14}
J.~Zhao and S.~Chen.
\newblock Guaranteed a posteriori error estimation for nonconforming finite
  element approximations to singularly perturbed reaction-diffusion problem.
\newblock {\em Adv. Comput. Math.}, 40:797--818, 2014.

\bibitem{ZS11}
L.~Zhu and D.~Sch\"otzau.
\newblock A robust a posteriori error estimate for hg-adaptive {D}{G} methods.
\newblock {\em IMA J. Num. Anal.}, 31:971--1005, 2011.

\end{thebibliography}

\def\cprime{$'$} \def\ocirc#1{\ifmmode\setbox0=\hbox{$#1$}\dimen0=\ht0
  \advance\dimen0 by1pt\rlap{\hbox to\wd0{\hss\raise\dimen0
  \hbox{\hskip.2em$\scriptscriptstyle\circ$}\hss}}#1\else {\accent"17 #1}\fi}
  \def\cprime{$'$}

\end{document}